\documentstyle[11pt]{article}

\catcode`\@=11

\newfam\msxfam
\newfam\msyfam

\def\hexnumber@#1{\ifcase#1 0\or1\or2\or3\or4\or5\or6\or7\or8\or9\or
	A\or B\or C\or D\or E\or F\fi }

\font\teneuf=eufm10
\font\seveneuf=eufm7
\font\fiveeuf=eufm5
\newfam\euffam
\textfont\euffam=\teneuf
\scriptfont\euffam=\seveneuf
\scriptscriptfont\euffam=\fiveeuf
\def\frak{\relaxnext@\ifmmode\let\next\frak@\else
 \def\next{\Err@{Use \string\frak\space only in math mode}}\fi\next}
\def\goth{\relaxnext@\ifmmode\let\next\frak@\else
 \def\next{\Err@{Use \string\goth\space only in math mode}}\fi\next}
\def\frak@#1{{\frak@@{#1}}}
\def\frak@@#1{\noaccents@\fam\euffam#1}

\edef\msx@{\hexnumber@\msxfam}
\edef\msy@{\hexnumber@\msyfam}

\mathchardef\boxdot="2\msx@00
\mathchardef\boxplus="2\msx@01
\mathchardef\boxtimes="2\msx@02
\mathchardef\square="0\msx@03
\mathchardef\blacksquare="0\msx@04
\mathchardef\centerdot="2\msx@05
\mathchardef\lozenge="0\msx@06
\mathchardef\blacklozenge="0\msx@07
\mathchardef\circlearrowright="3\msx@08
\mathchardef\circlearrowleft="3\msx@09
\mathchardef\rightleftharpoons="3\msx@0A
\mathchardef\leftrightharpoons="3\msx@0B
\mathchardef\boxminus="2\msx@0C
\mathchardef\Vdash="3\msx@0D
\mathchardef\Vvdash="3\msx@0E
\mathchardef\vDash="3\msx@0F
\mathchardef\twoheadrightarrow="3\msx@10
\mathchardef\twoheadleftarrow="3\msx@11
\mathchardef\leftleftarrows="3\msx@12
\mathchardef\rightrightarrows="3\msx@13
\mathchardef\upuparrows="3\msx@14
\mathchardef\downdownarrows="3\msx@15
\mathchardef\upharpoonright="3\msx@16

\mathchardef\downharpoonright="3\msx@17
\mathchardef\upharpoonleft="3\msx@18
\mathchardef\downharpoonleft="3\msx@19
\mathchardef\rightarrowtail="3\msx@1A
\mathchardef\leftarrowtail="3\msx@1B
\mathchardef\leftrightarrows="3\msx@1C
\mathchardef\rightleftarrows="3\msx@1D
\mathchardef\Lsh="3\msx@1E
\mathchardef\Rsh="3\msx@1F
\mathchardef\rightsquigarrow="3\msx@20
\mathchardef\leftrightsquigarrow="3\msx@21
\mathchardef\looparrowleft="3\msx@22
\mathchardef\looparrowright="3\msx@23
\mathchardef\circeq="3\msx@24
\mathchardef\succsim="3\msx@25
\mathchardef\gtrsim="3\msx@26
\mathchardef\gtrapprox="3\msx@27
\mathchardef\multimap="3\msx@28
\mathchardef\therefore="3\msx@29
\mathchardef\because="3\msx@2A
\mathchardef\doteqdot="3\msx@2B

\mathchardef\triangleq="3\msx@2C
\mathchardef\precsim="3\msx@2D
\mathchardef\lesssim="3\msx@2E
\mathchardef\lessapprox="3\msx@2F
\mathchardef\eqslantless="3\msx@30
\mathchardef\eqslantgtr="3\msx@31
\mathchardef\curlyeqprec="3\msx@32
\mathchardef\curlyeqsucc="3\msx@33
\mathchardef\preccurlyeq="3\msx@34
\mathchardef\leqq="3\msx@35
\mathchardef\leqslant="3\msx@36
\mathchardef\lessgtr="3\msx@37
\mathchardef\backprime="0\msx@38
\mathchardef\risingdotseq="3\msx@3A
\mathchardef\fallingdotseq="3\msx@3B
\mathchardef\succcurlyeq="3\msx@3C
\mathchardef\geqq="3\msx@3D
\mathchardef\geqslant="3\msx@3E
\mathchardef\gtrless="3\msx@3F
\mathchardef\sqsubset="3\msx@40
\mathchardef\sqsupset="3\msx@41
\mathchardef\vartriangleright="3\msx@42
\mathchardef\vartriangleleft="3\msx@43
\mathchardef\trianglerighteq="3\msx@44
\mathchardef\trianglelefteq="3\msx@45
\mathchardef\bigstar="0\msx@46
\mathchardef\between="3\msx@47
\mathchardef\blacktriangledown="0\msx@48
\mathchardef\blacktriangleright="3\msx@49
\mathchardef\blacktriangleleft="3\msx@4A
\mathchardef\vartriangle="0\msx@4D
\mathchardef\blacktriangle="0\msx@4E
\mathchardef\triangledown="0\msx@4F
\mathchardef\eqcirc="3\msx@50
\mathchardef\lesseqgtr="3\msx@51
\mathchardef\gtreqless="3\msx@52
\mathchardef\lesseqqgtr="3\msx@53
\mathchardef\gtreqqless="3\msx@54
\mathchardef\Rrightarrow="3\msx@56
\mathchardef\Lleftarrow="3\msx@57
\mathchardef\veebar="2\msx@59
\mathchardef\barwedge="2\msx@5A
\mathchardef\doublebarwedge="2\msx@5B
\mathchardef\angle="0\msx@5C
\mathchardef\measuredangle="0\msx@5D
\mathchardef\sphericalangle="0\msx@5E
\mathchardef\varpropto="3\msx@5F
\mathchardef\smallsmile="3\msx@60
\mathchardef\smallfrown="3\msx@61
\mathchardef\Subset="3\msx@62
\mathchardef\Supset="3\msx@63
\mathchardef\Cup="2\msx@64

\mathchardef\Cap="2\msx@65

\mathchardef\curlywedge="2\msx@66
\mathchardef\curlyvee="2\msx@67
\mathchardef\leftthreetimes="2\msx@68
\mathchardef\rightthreetimes="2\msx@69
\mathchardef\subseteqq="3\msx@6A
\mathchardef\supseteqq="3\msx@6B
\mathchardef\bumpeq="3\msx@6C
\mathchardef\Bumpeq="3\msx@6D
\mathchardef\lll="3\msx@6E

\mathchardef\ggg="3\msx@6F

\mathchardef\circledS="0\msx@73
\mathchardef\pitchfork="3\msx@74
\mathchardef\dotplus="2\msx@75
\mathchardef\backsim="3\msx@76
\mathchardef\backsimeq="3\msx@77
\mathchardef\complement="0\msx@7B
\mathchardef\intercal="2\msx@7C
\mathchardef\circledcirc="2\msx@7D
\mathchardef\circledast="2\msx@7E
\mathchardef\circleddash="2\msx@7F
\def\ulcorner{\delimiter"4\msx@70\msx@70 }
\def\urcorner{\delimiter"5\msx@71\msx@71 }
\def\llcorner{\delimiter"4\msx@78\msx@78 }
\def\lrcorner{\delimiter"5\msx@79\msx@79 }
\def\yen{\mathhexbox\msx@55 }
\def\checkmark{\mathhexbox\msx@58 }
\def\circledR{\mathhexbox\msx@72 }
\def\maltese{\mathhexbox\msx@7A }
\mathchardef\lvertneqq="3\msy@00
\mathchardef\gvertneqq="3\msy@01
\mathchardef\nleq="3\msy@02
\mathchardef\ngeq="3\msy@03
\mathchardef\nless="3\msy@04
\mathchardef\ngtr="3\msy@05
\mathchardef\nprec="3\msy@06
\mathchardef\nsucc="3\msy@07
\mathchardef\lneqq="3\msy@08
\mathchardef\gneqq="3\msy@09
\mathchardef\nleqslant="3\msy@0A
\mathchardef\ngeqslant="3\msy@0B
\mathchardef\lneq="3\msy@0C
\mathchardef\gneq="3\msy@0D
\mathchardef\npreceq="3\msy@0E
\mathchardef\nsucceq="3\msy@0F
\mathchardef\precnsim="3\msy@10
\mathchardef\succnsim="3\msy@11
\mathchardef\lnsim="3\msy@12
\mathchardef\gnsim="3\msy@13
\mathchardef\nleqq="3\msy@14
\mathchardef\ngeqq="3\msy@15
\mathchardef\precneqq="3\msy@16
\mathchardef\succneqq="3\msy@17
\mathchardef\precnapprox="3\msy@18
\mathchardef\succnapprox="3\msy@19
\mathchardef\lnapprox="3\msy@1A
\mathchardef\gnapprox="3\msy@1B
\mathchardef\nsim="3\msy@1C
\mathchardef\ncong="3\msy@1D

\mathchardef\varsubsetneq="3\msy@20
\mathchardef\varsupsetneq="3\msy@21
\mathchardef\nsubseteqq="3\msy@22
\mathchardef\nsupseteqq="3\msy@23
\mathchardef\subsetneqq="3\msy@24
\mathchardef\supsetneqq="3\msy@25
\mathchardef\varsubsetneqq="3\msy@26
\mathchardef\varsupsetneqq="3\msy@27
\mathchardef\subsetneq="3\msy@28
\mathchardef\supsetneq="3\msy@29
\mathchardef\nsubseteq="3\msy@2A
\mathchardef\nsupseteq="3\msy@2B
\mathchardef\nparallel="3\msy@2C
\mathchardef\nmid="3\msy@2D
\mathchardef\nshortmid="3\msy@2E
\mathchardef\nshortparallel="3\msy@2F
\mathchardef\nvdash="3\msy@30
\mathchardef\nVdash="3\msy@31
\mathchardef\nvDash="3\msy@32
\mathchardef\nVDash="3\msy@33
\mathchardef\ntrianglerighteq="3\msy@34
\mathchardef\ntrianglelefteq="3\msy@35
\mathchardef\ntriangleleft="3\msy@36
\mathchardef\ntriangleright="3\msy@37
\mathchardef\nleftarrow="3\msy@38
\mathchardef\nrightarrow="3\msy@39
\mathchardef\nLeftarrow="3\msy@3A
\mathchardef\nRightarrow="3\msy@3B
\mathchardef\nLeftrightarrow="3\msy@3C
\mathchardef\nleftrightarrow="3\msy@3D
\mathchardef\divideontimes="2\msy@3E
\mathchardef\varnothing="0\msy@3F
\mathchardef\nexists="0\msy@40
\mathchardef\mho="0\msy@66
\mathchardef\eth="0\msy@67
\mathchardef\eqsim="3\msy@68
\mathchardef\beth="0\msy@69
\mathchardef\gimel="0\msy@6A
\mathchardef\daleth="0\msy@6B
\mathchardef\lessdot="3\msy@6C
\mathchardef\gtrdot="3\msy@6D
\mathchardef\ltimes="2\msy@6E
\mathchardef\rtimes="2\msy@6F
\mathchardef\shortmid="3\msy@70
\mathchardef\shortparallel="3\msy@71
\mathchardef\smallsetminus="2\msy@72
\mathchardef\thicksim="3\msy@73
\mathchardef\thickapprox="3\msy@74
\mathchardef\approxeq="3\msy@75
\mathchardef\succapprox="3\msy@76
\mathchardef\precapprox="3\msy@77
\mathchardef\curvearrowleft="3\msy@78
\mathchardef\curvearrowright="3\msy@79
\mathchardef\digamma="0\msy@7A
\mathchardef\varkappa="0\msy@7B
\mathchardef\hslash="0\msy@7D
\mathchardef\hbar="0\msy@7E
\mathchardef\backepsilon="3\msy@7F
\def\Bbb{\relaxnext@\ifmmode\let\next\Bbb@\else
 \def\next{\Err@{Use \string\Bbb\space only in math mode}}\fi\next}
\def\Bbb@#1{{\Bbb@@{#1}}}
\def\Bbb@@#1{\noaccents@\fam\msyfam#1}

\catcode`\@=12

\font\tenfrak=eufm10
\font\sevenfrak=eufm7
\font\fivefrak=eufm5
\newfam\frakfam \def\frak{\fam\frakfam\tenfrak} \textfont\frakfam=\tenfrak
  \scriptfont\frakfam=\sevenfrak  \scriptscriptfont\frakfam=\fivefrak

\begin{document}
    \pagestyle{plain}
    \setlength{\baselineskip}{1.3\baselineskip}
    \setlength{\parindent}{\parindent}

\title{{\bf Hopf Algebroids and quantum groupoids}}
\author{Jiang-Hua Lu\\
Department of Mathematics, University of Arizona, Tucson 85721\\
jhlu@math.arizona.edu}
\maketitle

\newtheorem{thm}{Theorem}[section]
\newtheorem{lem}[thm]{Lemma}
\newtheorem{prop}[thm]{Proposition}
\newtheorem{cor}[thm]{Corollary}
\newtheorem{rem}[thm]{Remark}
\newtheorem{exam}[thm]{Example}
\newtheorem{nota}[thm]{Notation}
\newtheorem{dfn}[thm]{Definition}
\newtheorem{ques}[thm]{Question}
\newtheorem{eq}{thm}

\newcommand{\al}{\alpha}
\newcommand{\be}{\beta}
\newcommand{\rw}{\rightarrow}
\newcommand{\lrw}{\longrightarrow}
\newcommand{\Map}{\longmapsto}
\newcommand{\qed}{\begin{flushright} {\bf Q.E.D.}\ \ \ \ \
                  \end{flushright} }
\newcommand{\beqa}{\begin{eqnarray*}}
\newcommand{\eeqa}{\end{eqnarray*}}

\newcommand{\ot}{\mbox{$\otimes$}}
\newcommand{\otm}{\mbox{$\ot_{\scriptscriptstyle A}$}}
\newcommand{\xa}{\mbox{$x_{(1)}$}}
\newcommand{\xb}{\mbox{$x_{(2)}$}}
\newcommand{\xc}{\mbox{$x_{(3)}$}}
\newcommand{\ya}{\mbox{$y_{(1)}$}}
\newcommand{\yb}{\mbox{$y_{(2)}$}}
\renewcommand{\aa}{\mbox{$a_{(1)}$}}
\newcommand{\ab}{\mbox{$a_{(2)}$}}
\newcommand{\ac}{\mbox{$a_{(3)}$}}
\newcommand{\ba}{\mbox{$b_{(1)}$}}
\newcommand{\bt}{\mbox{$b_{(2)}$}}
\newcommand{\bc}{\mbox{$b_{(3)}$}}
\newcommand{\ca}{\mbox{$c_{(1)}$}}
\newcommand{\cb}{\mbox{$c_{(2)}$}}
\newcommand{\cc}{\mbox{$c_{(3)}$}}
\newcommand{\ha}{\mbox{$h_{(1)}$}}
\newcommand{\hb}{\mbox{$h_{(2)}$}}
\newcommand{\hc}{\mbox{$h_{(3)}$}}

\newcommand{\ts}{\mbox{$\tilde{\sigma}$}}
\newcommand{\las}{\mbox{${}_{\sigma}\!A$}}
\newcommand{\ras}{\mbox{$A_{\sigma}$}}
\newcommand{\rds}{\mbox{$\cdot_{\sigma}$}}
\newcommand{\lds}{\mbox{${}_{\sigma}\!\cdot$}}

\newcommand{\rhu}{\mbox{$\rightharpoonup$}}
\newcommand{\lhu}{\mbox{$\leftharpoonup$}}
\newcommand{\bb}{\mbox{$\bar{\beta}$}}
\newcommand{\bg}{\mbox{$\bar{\gamma}$}}
\newcommand{\la}{\mbox{$\langle$}}
\newcommand{\ra}{\mbox{$\rangle$}}

\newcommand{\id}{\mbox{${\em id}$}}
\newcommand{\Fun}{\mbox{${\em Fun}$}}
\newcommand{\End}{\mbox{${\em End}$}}
\newcommand{\Hom}{\mbox{${\em Hom}$}}
\newcommand{\ta}{\mbox{${\mbox{$\scriptscriptstyle A$}}$}}
\newcommand{\ap}{\mbox{$A_{\mbox{$\scriptscriptstyle P$}}$}}
\newcommand{\tx}{\mbox{$\mbox{$\scriptscriptstyle X$}$}}
\newcommand{\pp}{\mbox{$\pi_{\mbox{$\scriptscriptstyle P$}}$}}
\newcommand{\asemi}{\mbox{$\ap \#_{\sigma} A^*$}}
\newcommand{\dsemi}{\mbox{$A \#_{\Delta} A^*$}}

\newcommand{\fg}{\mbox{${\frak g}$}}
\newcommand{\uhg}{\mbox{$U_h {\frak g}$}}
\newcommand{\ugh}{\mbox{$U {\frak g} [[h]]$}}

\newcommand{\Gs}{\mbox{$G^*$}}
\newcommand{\pis}{\mbox{$\pi_{\sigma}$}}

\begin{abstract}
We introduce the notion of Hopf algebroids, in which neither
the total algebras nor the base algebras are
required to be commutative. We give a class of
Hopf algebroids associated to module algebras of the
Drinfeld doubles of Hopf algebras when the $R$-matrices
act properly. When this construction
is applied to quantum groups, we get examples of quantum groupoids,
which are semi-classical limits of Poisson groupoids. The example of
quantum $sl(2)$ is worked out in details.
\end{abstract}

\tableofcontents

\section{Introduction}
\label{sec_intro}

The notion of commutative Hopf algebroids was given in \cite{ra:cobor}.
It is modeled on the space of functions on a groupoid. In this
definition, both the algebra for the total space and the algebra
for the base space are commutative. The commutativity
of these algebras plays a crucial role in this definition.

It is natural to try to replace the algebras for both the total space
and the base space by non-commutative
algebras.
In \cite{ma:q-gpoid}, Maltsiniotis gives a definition
of Hopf algebroids, where the base algebras are still required to
be commutative and their images under the source and target maps
are required to lie in the centers of the total algebras.
See also \cite{b-m:quasi} for a study of
``quasi-quantum groupoids",
where the base algebras are again  commutative.

\bigskip
In this paper, we formulate a definition of
Hopf algebroids where none of the algebras are required to be commutative.
Like in \cite{ma:q-gpoid}, we use the term Hopf algebroids in general
and reserve
the term ``quantum groupoids" for the examples
associated to quantum groups.
A class of examples associated to module algebras
of the Drinfeld doubles of Hopf algebras is described in
Theorem \ref{thm_main}. An example
constructed from the quantum group $SL_q(2)$ when $q$ is a root of unity
is described in details.

\bigskip
Our work is motivated by the
the notion of {\bf Poisson groupoids} in Poisson geometry \cite{we:coi}.
Roughly speaking, we think of Hopf algebroids as ``quantizations" of
Poisson groupoids, just like quantum groups are quantizations of
Poisson groups. Our main example (see Theorem below) describes the
direct Hopf analog of some Poisson groupoids associated to
Poisson Lie groups.

Our definition reduces to the one given by Maltsiniotis in
\cite{ma:q-gpoid} when the base algebra is commutative and when
its images under the source and target
maps lie in the center of the total algebra. The semi-classical limits of
the quantum groupoids given by Maltsiniotis are Poisson groupoids
with the zero Poisson structure on the base spaces, while those
by our definition have arbitrary Poisson groupoids as
their semi-classical limits.

We would also like to point out that our definition of a
Hopf algebroid is not the same as the one given in
\cite{va:q-gpoid}, either, where Vainerman studies the Hopf analog of
transformation groupoids.
The semi-classical limit of his example does not give
rise to a Poisson groupoid. The main difference between  our example
and the one in \cite{va:q-gpoid} is that he uses
a direct product algebra structure while we use a smash
product algebra structure.

\bigskip
For any algebra $A$ over a field $k$, we can associate to it
the bi-algebroid $End_k(A)$ over $A$, and the
Hopf algebroid
$A \ot A^{op}$ over $A$. The latter is
obviously the Hopf analog of the coarse groupoid
structure on $X \times X$ over any space $X$.
Our main examples of Hopf algebroids,
however, are the following:

\bigskip
\noindent
{\bf Theorem}
{\it
Let $A$ be a Hopf algebra and let $D(A)$ be the Drinfeld double
of $A$. Let $V$ be a left $D(A)$-module algebra. Assume that the $R$-matrix
of $D(A)$ acts on $V \ot V$ in the following way:
\[
m^{op}_{\scriptscriptstyle V} \circ R ~ = ~ m_{\scriptscriptstyle V},
\]
where $m_{\scriptscriptstyle V}: V \ot V \rightarrow V$ is the
product on $V$ and $m^{op}_{\scriptscriptstyle V}$ is its opposite.
Then there is a Hopf algebroid structure over $V$
on the smash product algebra $H = V \# A$, where the smash product is
constructed using the left action of $A$ on $V$ that is the restriction
to $A$ of the action of $D(A)$ on $V$.
}

\bigskip
One example of such a $V$ is the dual Hopf algebra $A^*$ of $A$.
The smash product algebra $H = A^* \# A$ is, by definition, the
Heisenberg double of $A^*$, and we denote it by $H(A^*)$. Thus,
the above Theorem says that there
is a Hopf algebroid structure over $A^*$ on the Heisenberg double
$H(A^*)$.

\bigskip
In formulating the
definition, we use the following simple dictionary of Poisson geometry and
the theory of algebras (see \cite{lu:q-mom} for more details):

\bigskip
\begin{tabular}{ll} \vspace{2mm}
Poisson manifold $P$$~$& associative algebra $V$\\ \vspace{2mm}
Poisson morphisms$~$&algebra homomorphisms\\ \vspace{2mm}
Poisson submanifolds of $P$ $~$&two-sided ideas of $V$ \\ \vspace{2mm}
coisotropic submanifolds of $P$$~$&one-sided ideas of $V$ \\ \vspace{2mm}
Poisson group $G$ $~$&Hopf algebra $A$\\ \vspace{2mm}
Poisson action of $G$ on $P$ $~$&$A$-co-module algebra $V$\\ \vspace{2mm}
dual Poisson group $G^*$ $~$&dual Hopf algebra $A^*$\\ \vspace{2mm}
Semidirect product Poisson structure on $P \times G^*$ $~$ &smash
product on $V \ot A^*$\\ \vspace{2mm}
Poisson double $(D = G \bowtie G^*, \pi_{-})$ $~
$& dual of Drinfeld double $D(A)$ \\ \vspace{2mm}
symplectic double $(D = G \bowtie G^*, \pi_{+})$ $~$&Heisenberg double $H(A)$\\
\end{tabular}

\bigskip
We now recall the definition of a Poisson groupoid.
A {\bf groupoid} over a set $P$ is a set $X$ together with

1)  surjections $\alpha, \beta: X \longrightarrow P $ (called
the source and target maps respectively);

2)  $m: X_2 \longrightarrow X$ (called the multiplication), where
$X_2:= \{ (x, y) \in X \times X| \beta(x) = \alpha(y) \}$, such
that $\alpha(m(x,y)) = \alpha(x)$ and $\beta(m(x,y)) = \beta(y)$;

3) an injection $\epsilon: P \longrightarrow X$ (identities) such that
$\beta \epsilon = \alpha \epsilon = id_{\scriptscriptstyle P}$;

4) a bijection $\iota: X \longrightarrow X$ (inversion).

These maps must satisfy

 (1) associative law: $m (m (x, y), z) ~ = ~ m (x, m(y, z))$ (if one
side is defined, so is the other);

 (2) identities: for each $x \in X,  (\epsilon(\alpha (x)), x) \in X_{2},
(x, \epsilon (\beta(x))) \in X_{2}$ and $m (\epsilon(\alpha(x)), x) =
m(x, \epsilon (\beta(x))) = x$;

(3) inverses: for each $x \in X, (x, \iota(x)) \in X_{2},
(\iota (x), x) \in X_{2}, m(x, \iota(x)) = \epsilon (\alpha(x)),$
and $m(\iota (x), x) = \epsilon (\beta(x)).$

\bigskip
A groupoid $P$ over $X$ is called a Lie groupoid if both $X$ and $P$ are
smooth manifolds and if all the structure maps are smooth. We also require
that both $\alpha$ and $\beta$ be submersions. A {\bf Poisson groupoid}
\cite{we:coi}
is a Lie groupoid $P$ over $X$ together with a Poisson structure such that
the graph of the groupoid multiplication, i.e., the set
\[
\{~(m(x,y), ~ x, ~ y): ~~ \beta(x) = \alpha(y)~\} \subset X \times
X \times X
\]
is a coisotropic submanifold of $X \times {\bar X} \times {\bar
X}$,
where ${\bar X}$ denotes $X$ equipped with the negative Poisson
structure.  When the Poisson structure on $X$ is
nondegenerate, we say that $X$ is a {\bf symplectic groupoid}
\cite{we:sympoid} over $P$.

\bigskip
We summarize the properties of Poisson groupoids in the following
proposition. See \cite{we:coi} for the proofs. This proposition explains
some of the requirements we put in our definition of a Hopf algebroid.

\begin{prop}
\label{prop_poi}
Let $X$ be a Poisson groupoid over $P$ with structure maps
$\al, \be, m, \epsilon$ and $\iota$. Then

1) there is a unique Poisson structure on $P$ such that $\alpha$
is a Poisson map and $\beta$ is an anti-Poisson map. Consequently, the
submanifold
\[
\{(x, y): ~ \be(x) = \al(y) \} \subset X \times X
\]
is a coisotropic submanifold of $X \times X$;

2) the antipode map $\iota$ is an anti-Poisson isomorphism for $X$;

3) the map $\epsilon: P \rw X$ embedds $P$ into $X$ as a coisotropic
submanifold.
\end{prop}

Our main example of a Hopf algebroid is a direct Hopf analog of the
following theorem on Poisson Lie groups \cite{lu:thesis}.

\begin{thm}
\label{thm_p-gpoid}
Let $P$ be a Poisson manifold. Let $G$ be a Poisson Lie group, $G^*$
its dual Poisson Lie group and $(D, \pi_{-})$ the Poisson double
of $G$ and $G^*$. Assume that

1) there exists a right Poisson action $\sigma$ of $(D, \pi_{-})$ on
$P$:
\[
\sigma: ~ P \times D \longrightarrow P
\]

2) the Poisson structure on $P$ is determined by the action $\sigma$ in
the sense that the Poisson bracket of  two
functions $\phi_1$ and $\phi_2$ on $P$ is given by
\[
\{ \phi_1, ~ \phi_2 \} ~ = ~ (\sigma_{\phi_1}, ~ \sigma_{\phi_2}^{'})
\]
where $( ~ , ~ )$ denotes the pairing between the Lie algebra ${\frak g}$
of $G$ and the Lie algebra ${\frak g}^*$ of $G^*$, and
$\sigma_{\phi_1}$ and $\sigma_{\phi_2}^{'}$ are the ${\frak g}$
and ${\frak g}^*$-valued functions on $P$ respectively given by
\begin{eqnarray*}
& & \sigma_{\phi_1}(X)(p) ~ = ~ {\frac{d}{dt}}|_{t=0} \phi_1 (p \cdot \exp tX),
{}~~~ X \in {\frak g}\\
& & \sigma_{\phi_2}^{'}(\xi)(p) ~ = ~ {\frac{d}{dt}}|_{t=0} \phi_2 (p \cdot
\exp
t\xi), ~~~ \xi \in {\frak g}^*.
\end{eqnarray*}
Then the manifold $P \times G^*$, together with the transformation
groupoid structure over $P$ defined by the restricted $G^*$
action on $P$ and the semi-Poisson structure  defined by the
restricted action of $G$ on $P$, becomes a Poisson groupoid
over $P$.
\end{thm}

Throughout this paper, vector spaces are
assumed to be over the ground field $k$, and tensor products
are over $k$ unless otherwise indicated. Algebras are assumed to
have identity elements, and algebra morphisms are assumed to
preserve identity elements.

\bigskip
The paper is organized as follows. In Section \ref{sec_bi-agoid},
we give the definition of a bi-algebroid. We show
in Section \ref{sec_two-exam} that for any
finite dimensional algebra $A$, the space $End_k(A)$ has a
natural bi-algebroid structure over $A$. Moreover, for any
other bi-algebroid $H$ over $A$, there is a bi-algebroid morphism
from $H$ to $End_k(A)$. Another example of a bi-algebroid
over $A$ is $A \otimes A^{op}$. In Section \ref{sec_hopf-agoid},
we give the definition of a Hopf algebroid. The main examples are
described in Theorem \ref{thm_main} in Section \ref{sec_main}.
The Hopf algebroid structure
on the Heisenberg double $H(A^*)$ of a Hopf algebra $A^*$ is
described in details in Section \ref{sec_hsberg}. Finally,
we give the example related to the quantum group $SL_q(2)$
at a root of unity in Section \ref{sec_q-poid}.

\bigskip
\noindent
{\bf Acknowledgement} The author would like to thank Alan Weinsterin,
Gorges Maltsiniotis, Alain Bruguieres and Yvette Kosmann-Schwatzbach
for helpul discussions. Part of the work was done at the Institute
for Advanced Study at Princeton, the Emile Borel Center in Paris and
the Isaac Newton Institute for Mathematical Sciences in
Cambridge, England, to which the author
expresses her gratitude for their hospitality.

\section{Bi-algebroids}
\label{sec_bi-agoid}

\begin{dfn}
\label{dfn_algebroid}
{\em
A  bi-algebroid consists of the following data:

1) an algebra $H$ called the {\bf total algebra};

2) an algebra $A$ called the {\bf base algebra};

3) the {\bf source map}: an algebra homomorphism
\[
\al: ~ A \lrw H
\]
and the {\bf target map}: an algebra anti-homomorphism
\[
\be: ~ A \lrw H
\]
such that the images of $\al$ and $\be$ commute in $H$, i.e., $\forall
a, b \in A$,
\begin{equation}
\label{eq_albe}
\al(a) \be(b) ~ = ~ \be(b) \al (a).
\end{equation}
There is then a natural $(A, A)$-bimodule structure on $H$ given by
\begin{eqnarray}
\label{eq_lambda}
& \lambda: & A \ot H \lrw H: ~~ a \ot h \Map \al (a)h\\
\label{eq_rho}
& \rho: & H \ot A \lrw H: ~~ h \ot a \Map \be(a) h.
\end{eqnarray}
Using this bimodule structure, we can form the $(A, A)$-bimodule
product $H \ot_{\scriptscriptstyle A} H$ of $H$ with itself. As
a vector space over $k$, it is simply the quotient of $H \ot H$ by
the subspace
\begin{equation}
\label{eq_ir}
I_2 := \left \{ \be(a) h_1 \ot h_2  -  h_1 \ot \al(a) h_2 = (\be(a) \ot 1
- 1 \ot \al(a) )(h_1 \ot h_2), ~~ h_1, h_2 \in H, a \in A \right \}.
\end{equation}
Notice that in this case,  $I_2$ is in fact the right ideal of $H \ot H$
generated by the subset
\[
\{\be(a) \ot 1 - 1 \ot \al(a): ~ a \in A \}.
\]
We will still denote elements of $H \ot_{\scriptscriptstyle A} H$ by
$h_1 \ot h_2$. The $(A, A)$-bimodule structure on
$H \ot_{\scriptscriptstyle A} H$ is now given by
\[
A \ot (H \ot_{\scriptscriptstyle A} H) \lrw H \ot_{\scriptscriptstyle A} H:
{}~~ a \ot (h_1 \ot h_2) \Map \al(a)h_1 \ot h_2
\]
and
\[
(H \ot_{\scriptscriptstyle A} H) \ot A \lrw H \ot_{\scriptscriptstyle A} H: ~~
(h_1 \ot h_2) \ot a \Map h_1 \ot \be(a) h_2.
\]
We can then form the $(A, A)$-bimodule product of $H \ot_{\scriptscriptstyle A}
H$ with $H$ to get the triple product
$H \ot_{\scriptscriptstyle A} H \ot_{\scriptscriptstyle A} H$.
In general, we can form the $n$'th $(A, A)$-bimodule
product of $H$ with itself. It is the quotient of $H^{\otimes n}$ by
the right ideal $I_n$ generated
by the subspace
\begin{equation}
\label{eq_In}
\{\underbrace{1 \ot \cdots \ot 1}_{i-1} \ot \be(a_i) \ot 1 \ot \cdots \ot 1 -
\underbrace{ 1 \ot \cdots \ot 1}_{i-1} \ot 1 \ot  \al(a_i) \ot 1 \ot
\cdots \ot 1, ~ i = 1, ..., n-1, a_i \in A\}.
\end{equation}

4) the {\bf co-product}: an $(A, A)$-bimodule map
\[
\Delta: ~ H \lrw H \ot_{\scriptscriptstyle A} H: ~~ h \Map \ha \ot \hb
\]
with $\Delta(1) = 1 \ot 1$ satisfying the co-associativity:
\[
(\Delta \ot_{\scriptscriptstyle A} id_{\scriptscriptstyle H} ) \Delta
{}~ = ~ (id_{\scriptscriptstyle H} \ot_{\scriptscriptstyle A} \Delta) \Delta: ~
H \longrightarrow H \otm H \otm H.
\]
The co-product $\Delta$ and the algebra structure on $H$ are
required to be
{\bf compatible} in the sense that the kernel of the following map
\begin{equation}
\label{eq_phi}
\Phi:~
H \ot H \ot H \lrw H \ot_{\scriptscriptstyle A} H: ~~ h_1 \ot h_2 \ot h_3
\Map (\Delta h_1) (h_2 \ot h_3)
\end{equation}
is a {\bf left ideal} of $H \ot H^{op} \ot H^{op}$, where $H^{op}$
denotes $H$ with the opposite product. Here we are using the
fact that $H \ot H$ acts on $H \ot_{\scriptscriptstyle A} H$ from
the right by right multiplications;

5) the co-unit map: an (A, A)-bimodule map
\[
\epsilon: ~~ H \lrw A
\]
satisfying
\[
\epsilon(1_{\scriptscriptstyle H}) ~ = ~ 1
\]
(it follows then that $\epsilon \be  =  \epsilon \al  =  id_{\scriptscriptstyle
A}$) and
\begin{equation}
\label{eq_co-unit}
\lambda (\epsilon \ot id_{\scriptscriptstyle H}) \Delta ~ = ~ \rho
(id_{\scriptscriptstyle H} \ot \epsilon) \Delta ~ =
{}~ id_{\scriptscriptstyle H}: ~ H \lrw H,
\end{equation}
where $\lambda$ and $\rho$ are respectively the left and right
$A$-module maps for $H$ given by (\ref{eq_lambda})
and (\ref{eq_rho}).
Notice that two maps on the left hand side of the above equation are
well-defined. It is also required to be {\bf compatible} with the algebra
structure on $H$ in the sense that the kernel of $\epsilon$ is a
{\bf left ideal} of $H$.
}
\end{dfn}

\bigskip
The following Lemma shows that when $A$ is the $1$-dimensional
algebra $k$, our definition of a bi-algebroid over $k$ is reduced to  that
of a bi-algebra over $k$.

\begin{lem}
\label{lem_morphi}
Let $A$ and $B$ be two algebras. A linear map $f: A \rw B$
mapping $1_{\scriptscriptstyle A}$ to $1_{\scriptscriptstyle B}$
is an algebra homomorphism if and only if the kernel of the
map
\[
F: ~ A \ot B \lrw B: ~ a \ot b \Map f(a)b
\]
is a left ideal of $A \ot B^{op}$. In this case, the left ideal is
generated by the subspace $\{a\ot 1 - 1 \ot f(a): a \in A\}$.
\end{lem}

\noindent
{\bf Proof.}
Assume first that $f$ is an algebra morphism. Then it is easy to see that
$\ker F$ is a left ideal of $A \ot B^{op}$.
Moreover, if
$\sum a_i \ot b_i \in \ker F$, then
\[
\sum a_i \ot b_i ~ = ~ \sum (1 \ot b_i) (a_i \ot 1 - 1 \ot f(a_i)) \in A \ot
B^{op}.
\]
Thus $\ker F$ coincides with the left ideal generated by the given subspace.
Conversely, assume that $\ker F$ is a left ideal of $A \ot B^{op}$.
Then for any $a_1, a_2 \in A$, since
\[
a_2 \ot 1 ~ - ~ 1 \ot f(a_2) \in \ker F,
\]
we know that
\[
a_1 a_2 \ot 1 - a_1 \ot
f(a_2) ~ = ~ (a_1 \ot 1) (a_2 \ot 1 - 1 \ot f(a_2))
\in \ker F.
\]
Thus $f(a_1 a_2) = f(a_1) f(a_2)$. Hence $f$ is an algebra
homomorphism.
\qed

\bigskip
Assume that $H$ is a bi-algebroid over the algebra $A$ with
structure maps $\al, \be, \Delta$ and $\epsilon$.
Then since  $\epsilon \al = id_{\scriptscriptstyle A}$, the map $\epsilon$
is subjective, so we can identify $A$ with the quotient space
$H / \ker \epsilon$. Since $\ker \epsilon \subset H$ is a left ideal
of $H$, there is an induced left action of $H$ on $A$.
Since $\alpha: A \lrw H$ is a section of $\epsilon$, we can identify $A$ with
the subspace $\al(A) \subset H$. The left
action of $H$ on $A$ can then be explicitly written down as
\begin{equation}
\label{eq_asso}
T_1: ~ H \lrw End_k(A):  ~ h(a) := \epsilon(h \al(a)).
\end{equation}

\begin{dfn}
\label{dfn_induced}
{\em
We call the action $T_1$ given by (\ref{eq_asso}) the left action of $H$
on $A$ associated to the bi-algebroid structure on $H$ over $A$.
}
\end{dfn}

\begin{prop}
\label{prop_asso}
The action $T_1$ given by (\ref{eq_asso}) has the following properties:

1) The composition $T_1 \circ \al$ (resp. $T_1 \circ \beta$): $A \rw
End_k(A)$ gives the left (resp. right) action of $A$ on itself
by left (resp. right) multiplications;

2) $\epsilon (h) ~ = ~ h(1_{\scriptscriptstyle A}), ~~~~\forall h \in H$;

3) The map $T_2: ~H \ot H \lrw Hom_k (A \ot A, A)$ given by
\[
T_2 (h_1 \ot h_2) (a \ot b) ~ = ~ h_1 (a) h_2 (b)
\]
induces a well-defined map, still denoted by $T_2$, from
$H \otm H$ to $Hom_k (A \ot A, A)$. We have
\begin{equation}
\label{eq_comp}
h(ab) ~ = ~ (T_2 \circ \Delta) (h) (a \ot b) ~ = ~
\ha(a) ~ \hb(b), ~~~ \forall h \in H, ~ a, b \in A.
\end{equation}
\end{prop}

\noindent
{\bf Proof.} Statement 1) follows from the
fact that $\epsilon: H \rw A$ is a $(A, A)$-bimodule map, i.e.,
for all $a \in A$ and $h \in H$,
\[
\epsilon(\al(a) h) ~ = ~ a \epsilon(h), ~~~ \epsilon(\be(a)h) ~ = ~ \epsilon
(h) a.
\]
Statement 2) follows from the definition of the action $T_1$. Since
$T_2$ maps elements in the right ideal $I_2$ given by (\ref{eq_ir})
to zero, it induces a well-defined map from $H \otm H$ to $Hom_k(A \ot A, A)$.
To see the second part of statement 3), we first prove the following:
\begin{equation}
\label{eq_change1}
h \al(a) ~ = ~ \al (\ha(a)) \hb, ~~~~ \forall a \in A, h \in H.
\end{equation}
Since $\Delta(1) = 1 \ot 1$ and since $\Delta$ is a left $A$-module map,
we have
\[
\Delta(\al(a)) ~ = ~ \al(a) \ot 1, ~~~ a \in A.
\]
Thus
\[
\al(a) \ot 1 \ot 1 - 1 \ot \al(a) \ot 1 \in \ker \Phi,
\]
where $\Phi$ is given by (\ref{eq_phi}). Since $\ker \Phi$ is
a left ideal in $H \ot {\bar H} \ot {\bar H}$, we have
\[
h \al(a) \ot 1 \ot 1 - h \ot \al(a) \ot 1 \in \ker \Phi, ~~~ h \in H.
\]
Hence
\[
\Delta(h\al(a)) ~ = ~ \ha \al(a) \ot \hb.
\]
Applying (\ref{eq_co-unit}), we get
\[
\alpha(\epsilon(\ha \alpha(a))) \hb ~ = ~ h \alpha(a).
\]
But
\[
\epsilon (\ha \alpha(a)) ~ = ~ \ha(a)
\]
by definition. Therefore we get
(\ref{eq_change1}).
Consequently, we have
\[
h(ab) ~ = ~ (h \al(a))(b) ~ = ~ \ha(a) \hb(b).
\]
\qed

\begin{rem}
\label{rem_conditions}
{\em
Note how the condition that $\ker \Phi$ be a left idea in $H \ot \bar{H}
\ot \bar{H}$ and the compatibility condition of $\epsilon$ and $\Delta$
are used in the above proof.
}
\end{rem}

\begin{rem}
\label{rem_what2}
{\em
Apart from the right $A$-module structure $\rho$ on $H$ given by
(\ref{eq_rho}), there is another right $A$-module structure on $H$
given by
\[
H \ot A \lrw H: ~ h \ot a \Map h \al(a).
\]
Equation (\ref{eq_change1}) should be thought of as expressing this
right $A$-module structure in terms of the left $A$-module structure
$\lambda$ as defined in (\ref{eq_lambda}).
Similarly, we have
\begin{equation}
\label{eq_change2}
h \be(a) ~ = ~ \be(\hb(a)) \ha, ~~~ \forall a \in A, ~ h \in H,
\end{equation}
which expresses the left $A$ module structure on $H$ given by
\[
A \ot H \lrw H: ~ a \ot h \Map h \be(a)
\]
in terms of the right $A$-module map $\rho$.
}
\end{rem}

We think of (\ref{eq_comp}) as the ``product" rule for the action
$T_1$ of $H$ on $A$ with respect to the algebra structure of $A$.
Thus, roughly speaking, a bi-algebroid over $A$ is an algebra
which acts on $A$ in such a way that
1) it contains as part of it the left and right
actions of $A$ on itself by left and right multiplications, and 2) it
obeys a product rule specified by the co-product of the bi-algebroid.

\section{Two examples of bi-algebroids: $A \ot A^{op}$ and $End_k(A)$}
\label{sec_two-exam}

We now give two examples of bi-algebroids.

\begin{exam}
\label{exam_coarse}
{\em
Let $A$ be an arbitrary algebra and let $H = A \ot A^{op}$, i.e.,
$H$ is the direct product algebra of $A$ and its opposite. Then
there is a bi-algebroid structure on $H$ over $A$ with the
following structure maps:

1) the source map $\al: A \rw H: ~ a \mapsto a \ot 1$, and the target map
$\be: A \rw H: a \mapsto 1 \ot a$;

2) the co-product
\[
\Delta: H \lrw H \otm H: ~~a \ot b \Map (a \ot 1) \ot (1 \ot b).
\]

3) the co-unit
\[
\epsilon: ~ H \lrw A: ~ a \ot b \Map ab.
\]

It is easy to check that these maps indeed define a Hopf algebroid
structure on $H = A \ot A^{op}$ over $A$. The associated left action
of $H$ on $A$ is given by
\[
T_1: ~ H \lrw End_k (A): ~ T_1 (a \ot b) (c) ~ = ~ acb.
\]
}
\end{exam}

\begin{rem}
\label{rem_coarse}
{\em
This example is modeled on the Poisson groupoid structure on
$P \times {\bar P}$ for any Poisson manifold (see \cite{we:coi}),
which is in turn a Poisson version of the coarse groupoid
structure on $X \times X$ over any space $X$. Thus we
call $A \ot A^{op}$ the
``coarse Hopf algebroid" of $A$.}
\end{rem}

\begin{rem}
\label{rem_coarse2}
{\em
We will show in Section \ref{sec_hopf-agoid} that
in addition to the bi-algebroid
structure, there is actually a Hopf algebroid structure
on $A \ot A^{op}$.
}
\end{rem}

\begin{exam}
\label{exam_end}
{\em
Let $A$ be any finite dimensional algebra over $k$. Let $H = End_k(A)$
be the algebra of $k$-linear maps from $A$ to itself. We will show that
there is a bi-algebroid structure on $H$ over $A$. We need the
following Lemma.

\begin{lem}
\label{lem_end}
Let $A$ and $H$ be as above.
For any $a \in A$, let $\al(a) \in H$ and $\be(a) \in H$ be respectively the
maps from $A$ to itself defined
by the left and right multiplications by $a$. Then the
right ideal $I_2$ of $H \ot H$ given by (\ref{eq_ir})
coincides with the kernel of the map
\[
T_2: ~ H \ot H \lrw Hom_k(A \ot A, A): ~ T_2(h_1 \ot h_2) (a \ot b) = h_1(a)
h_2(b).
\]
Thus we get a natural indentification between $H \otm H$ and the
space $Hom_k (A \ot A, A)$. In general, for any positive integer
$n$, the map
\[
T_n: ~ H^{\otimes n} \lrw Hom_k(A^{\otimes n}, ~ A): ~
T_n(h_1 \ot \cdots \ot h_n) (a_1 \ot \cdots \ot a_n) ~ = ~
h_1(a_1) \cdots h_n(a_n)
\]
gives a natural identification of the
$(A, A)$-bimodules
$\underbrace{H \otm \cdots \otm H}_n$ and
$Hom_k(A^{\otimes n}, ~A)$.
\end{lem}

\noindent
{\bf Proof.} Clearly, $I_2 \subset \ker T_2$. Conversely, let $\{a_s\}$ be a
basis of
$A$ and let $\{x_s\}$ be its dual basis for $A^*$.
For each $s = 1, 2, ..., \dim A$, let
$h_s \in End_k(A)$ be defined by
\[
h_s (a) ~ = ~ (a, ~ x_s) 1_{\scriptscriptstyle A}.
\]
Then any $\sum_i h_{i}^{'} \ot h_{i}^{''} \in \ker T_2$ can be written as
\[
\sum_i h_{i}^{'} \ot h_{i}^{''} ~ = ~ \sum_i \sum_{s = 1}^{\dim A}
\left( \be(h_{i}^{'} (a_s)) \ot 1 - 1 \ot \al(h_{i}^{'}(a_s)) \right) (h_s \ot
h_{i}^{''}).
\]
Thus $\sum_i h_{i}^{'} \ot h_{i}^{''} \in I_2$. Hence $I_2 = \ker T_2$. For a
general
$n \geq 2$, recall that the space $H \otm H \otm \cdots \otm H$ is the quotient
of $H^{\otimes n}$ by the right ideal $I_n$ given in (\ref{eq_In}).
Clearly, $I_n \subset \ker T_n$. Conversely, any $\sum_i h_{i}^{'} \ot
h_{i}^{''}
\ot \cdots \ot h_{i}^{(n)} \in \ker T_n$ can be written as
\[
\sum_i h_{i}^{'} \ot h_{i}^{''}
\ot \cdots \ot h_{i}^{(n)} ~ = ~ \sum_i \sum_{j=1}^{n-1}
\sum_{1 \leq s_1, \cdots, s_j \leq \dim A}
(h_{j;~ s_1, \cdots, s_j}) (h_{s_1} \ot \cdots \ot h_{s_j}
\ot h_{i}^{(j+1)} \ot \cdots \ot h_{i}^{(n)}),
\]
where for $1 \leq s_1, \cdots, s_j \leq \dim A$,
\beqa
 h_{j; ~s_1, \cdots, s_j} & = &
 \underbrace{1 \ot \cdots \ot 1}_{j-1}
 \ot \left(\be(h_{i}^{'}(a_{s_1})
h_{i}^{''}(a_{s_2}) \cdots h_{i}^{(j)} (a_{s_j})) \ot 1\right)
\ot 1 \ot \cdots \ot 1 \\
& &~- ~ \underbrace{1 \ot \cdots \ot 1}_{j-1}\ot \left( 1 \ot
\al(h_{i}^{'}(a_{s_1})
h_{i}^{''}(a_{s_2}) \cdots h_{i}^{(j)} (a_{s_j}))\right) \ot 1 \ot \cdots \ot
1\\
& &~\in ~  I_n.
\eeqa
Thus
\[
\sum_i h_{i}^{'} \ot h_{i}^{''}
\ot \cdots \ot h_{i}^{(n)} \in I_n.
\]
Hence $I_n = \ker T_n$.
Consequently, the space $H \otm H \otm \cdots \otm H$
can be identified with the space $Hom_k(A^{\otimes n}, ~A)$.
Since $T_n$ is an $(A, A)$-bimodule map, this is an identification
of $(A, A)$-bimodules.
\qed

\bigskip
We now define the bi-algebroid structure on $End_k(A)$ over $A$.
Define
\[
\Delta: ~ H \lrw H \otm H \cong Hom_k (A \ot A, A)
\]
by
\[
\Delta(h) (a \ot b) ~ = ~ h(ab),
\]
and
\[
\epsilon:~ H \lrw A
\]
by
\[
\epsilon(h) ~= ~ h(1_{\scriptscriptstyle A}).
\]
Then it is straightforward to check that
the maps $\al, \be, \Delta$ and $\epsilon$ define a bi-algebroid
structure on $H = End_k (A)$ over $A$. Here we only prove the
compatibility of the product and the co-product on $H =
End_k (A)$. Identify
\[
H \otm H ~ \cong ~ Hom_k (A \ot A, ~ A)
\]
by the map $T_2$. The map $\Phi$ given in
(\ref{eq_phi}) now becomes
\[
\Phi: ~ H \ot H \ot H \lrw Hom_k(A \ot A, ~ A): ~
\Phi(h_1 \ot h_2 \ot h_3) (a \ot b) ~ = ~ h_1 \left( h_2(a) h_3(b)
\right).
\]
Let $\{a_s\}$ be a basis of $A$, and let $\{x_s\}$ be the
dual basis for $A^*$. Then the algebra
structure on $A$ corresponds to the element
\[
m ~ = ~ a_s a_t ~ \ot~  x_s ~ \ot ~x_t \in A \ot A^* \ot A^*
\]
under the natural identification of $A \ot A^* \ot A^*$ with
$Hom_k (A \ot A, A)$.
The kernel of $\Phi$ is now the annihilator in
$H \ot H^{op} \ot H^{op}$ of the element
$m \in A \ot A^* \ot A^*$ with respect to the left representation
of $H \ot H^{op} \ot H^{op}$ on $A \ot A^* \ot A^*$,
where $H^{op}$ acts on $A^*$
by dualizing the action of $H = End_k(A)$ on $A$. Hence
$ker \Phi$ is a left ideal in $H \ot H^{op} \ot H^{op}$.
}
\end{exam}

\begin{dfn}
\label{dfn_homo}
{\em
A bi-algebroid morphism from a bi-algebroid $(H_1, A_1, \al_1, \be_1,
\Delta_1, \epsilon_1)$ to another bi-algebroid
$(H_2, A_2, \al_2, \be_2, \Delta_2, \epsilon_2)$ consists of an
algebra morphism $T: H_1 \rw H_2$ and an algebra morphism
$t: A_1 \rw A_2$ witch commute with all the structure maps.
}
\end{dfn}

Our previous discussion leads immediately to the following proposition:

\begin{prop}
\label{prop_univ}
Assume that $A$ is a finite dimensional algebra. Then for any bi-algebroid
$H$ over $A$, there is a bi-algebroid morphism from $H$ to $End_k (A)$
given by the map $T_1$ defined by (\ref{eq_asso}) and the identity map
on $A$.
\end{prop}

\begin{exam}
\label{exam_simple}
{\em
When $A$ is a finite-dimensional simple algebra,
the map
\[
T_1: ~ A \ot A^{op} \lrw End_k(A): ~ T_1(a \ot b)(c) ~ = ~ acb
\]
defines an algebra isomorphism. Thus in this case, the two
bi-algebroids $A \ot A^{op}$ and $End_k(A)$ over $A$ are isomorphic.
}
\end{exam}

\section{Hopf algebroids}
\label{sec_hopf-agoid}

Naturally, Hopf algebroids should be bi-algebroids with antipodes.

\begin{dfn}
\label{dfn_hopf}
{\em
A Hopf algebroid is a bi-algebroid $H$ over an algebra $V$ with
structure maps $\al, \be, \Delta$ and $\epsilon$ together with a
bijective map $\tau: H \rw H$, called the antipode map,
which has the following properties:

1) $\tau$ is an algebra anti-isomorphism for $H$;

2) $\tau \be = \al$;

3) $m_{\scriptscriptstyle H} (\tau \ot id) \Delta ~ = ~
\be \epsilon \tau: ~ H \rightarrow H$,
where $m_{\scriptscriptstyle H}$ denotes
the multiplication of $H$;

4) there exists a linear map $\gamma: H \otm H \rw H \ot H$
with the following properties:

4a) $~\gamma$ is a
section for the natural projection $p: H \ot H \rw H \otm H$;

4b) the following identity holds:
\begin{equation}
\label{eq_section}
m_{\scriptscriptstyle H} (id \ot \tau ) \gamma \Delta  ~ = ~ \al \epsilon:
{}~ H \lrw H.
\end{equation}
}
\end{dfn}

\begin{rem}
\label{rem_section}
{\em
1)  The map $m_{\scriptscriptstyle H} (\tau \ot id) \Delta$ in the formula
in 3) is well-defined as a map from $H$ to $H$, but the map
$m_{\scriptscriptstyle H}(id \ot \tau) \Delta$ is not well-defined.
This is why we need to require the existence of the section $\gamma$.

2) In general, $\epsilon \tau \neq \epsilon$, as can be seen in
Example \ref{exam_hopf-coarse}.

3) In \cite{ma:q-gpoid}, Maltsiniotis studies the case
when $A$ is commutative and when the
images of $V$ under $\alpha$ and $\beta$ lie in the
center of $H$. In this case, his definition coincides with ours.
}
\end{rem}

\begin{prop}
\label{prop_d}
There exists an algebra automorphism $\theta$ for $V$ such that
\begin{equation}
\label{eq_d-exist}
\tau \al ~ = ~ \be \theta.
\end{equation}
\end{prop}

\noindent
{\bf Proof} ~ We know that
\[
\Delta \al(v) ~ = ~ \al(v) ~ \ot ~ 1.
\]
Apply the identity $m_{\scriptscriptstyle H} (\tau \ot id)
\Delta = \be \epsilon \tau$ to $\al (v)$. We get
\[
\tau (\al(v)) ~ = ~ \be \epsilon \tau \al (v).
\]
Set
\begin{equation}
\label{eq_d-dfn}
\theta(v) ~ = ~ \epsilon \tau \al (v) \in V.
\end{equation}
Then
\[
\tau (\al(v)) ~ = ~ \be (\theta(v)).
\]
Since $\tau$ is an algebra anti-isomorphism for $H$ and since $\be$
is injective, we know that $\theta$ is an algebra automorphism for $V$.
\qed

\begin{exam}
\label{exam_hopf-coarse}
{\em
For any algebra $A$, let $H = A \ot A^{op}$ be the bi-algebroid over
$A$ as described in Example \ref{exam_coarse}. Define $\tau:  H \rw H$ by
\[
\tau(a \ot b) ~ = ~ b \ot a.
\]
Then $H = A \ot A^{op}$ is a Hopf algebroid over $A$ with antipode $\tau$.
The induced algebra isomorphism $\theta$ for $A$ in this example is the
identity map.
}
\end{exam}

\section{The Hopf algebroid $V \# A$}
\label{sec_main}

We now give a class of Hopf algebroids. They should be
considered as the Hopf version of transformation groupoids.

\begin{thm}
\label{thm_main}
Let $A$ be a Hopf algebra and let $D(A)$ be the Drinfeld double
of $A$. Let $V$ be a left $D(A)$-module algebra. Assume that the $R$-matrix
of $D(A)$ acts on $V \ot V$ in the following way:
\begin{equation}
\label{eq_R}
m^{op}_{\scriptscriptstyle V} \circ R ~ = ~ m_{\scriptscriptstyle V},
\end{equation}
where $m_{\scriptscriptstyle V}: V \ot V \rightarrow V$ is the
product on $V$ and $m^{op}_{\scriptscriptstyle V}$ is its opposite.
Then there is a Hopf algebroid structure over $V$
on the smash product algebra $H = V \# A$, where the smash product is
constructed using the left action of $A$ on $V$ that is the restriction
to $A$ of the action of $D(A)$ on $V$.
\end{thm}

The rest of this section is devoted to the explanation and the proof of this
theorem.

\bigskip
We first recall the definition of the Drinfeld double $D(A)$ of
a Hopf algebra $A$. As  a vector space, $D(A)$ is isomorphic to the
tensor product space $A^* \ot A$. As a coalgebra, it has the direct
product coalgebra structure of $A^{*coop}$ and $A$, where
$A^{*coop}$ is $A^*$ with the opposite co-product. Its algebra
structure is defined as follows: for $x \ot a$ and $y \ot b$ in
$D(A)$, their product in $D(A)$ is given by
\[
D(A) \ni (x \ot a) (y \ot b) ~ = ~ x ~ (\aa \triangleright \yb ) ~ \ot ~ (\ab
\triangleleft \ya ) ~ b,
\]
where
\[
A \ot A^* \longrightarrow A^*: ~ a \ot x \longmapsto
a \triangleright x ~ = ~ \aa  \rhu ~x ~\lhu S^{-1} (\ab)
\]
is the left co-adjoint representation of $A$ on $A^*$ and
\[
A \ot A^* \longrightarrow A: ~ a \ot x \longmapsto
a \triangleleft x ~ = ~ S^{-1}(\xa) \rhu ~a ~\lhu \xb
\]
is the right co-adjoint representation of $A^*$ on $A$. Here we
are using the standard notion $\rhu$ and $\lhu$
to denote the left and right representations of $A$ and $A^*$
on each other given as follows: for $a \in A$ and $x \in A^*$,
\begin{eqnarray}
\label{eq_rhu}
& & a ~\rhu ~ x ~ = ~ \xa \la\xb, ~ a\ra\\
\label{eq_lhu}
& &x ~\lhu ~a ~ = ~ \xb \la a, ~ \xa \ra.
\end{eqnarray}
The antipode map $S_{\scriptscriptstyle D(A)}$ of $D(A)$
is given by
\[
S_{\scriptscriptstyle D(A)} (x \ot a) ~ = ~ S(a) S^{-1}(x).
\]
Notice that the natural inclusions
\beqa
& &A \hookrightarrow D(A): ~~ a \Map 1 \ot a\\
& &A^* \hookrightarrow D(A): ~~ x \Map x \ot 1
\eeqa
are Hopf algebra morphisms.

The $R$-matrix for $D(A)$ is given by
\begin{equation}
\label{eq_R-dfn}
R ~ = ~ (1 \ot a_t) \ot (x_t \ot 1) \in D(A) \ot D(A)
\end{equation}
where $\{a_t\}$ is a basis for $A$, and $\{x_t\}$ is its dual
basis for $A^*$. Condition (\ref{eq_R}) now reads as
\begin{equation}
\label{eq_R1}
x_t(u) ~ a_t(v) ~ = ~ vu, ~~~~~\forall u, v \in V.
\end{equation}

We also recall that $V$ is said to be a left $D(A)$-module algebra if the
action of $D(A)$ on $V$ satisfies
\[
d(uv) ~ = ~ d_{(1)} (u) ~ d_{(2)} (v), ~~~~~\forall d \in D(A), ~ u, v \in V.
\]
In particular, since the inclusion
\[
A \lrw D(A): ~ a \Map 1 \ot a
\]
is a Hopf algebra inclusion, we can restrict the action to an action
of $A$ on $V$, making $V$ into a left $A$-module algebra. The smash
product algebra $H = V \# A$ is, by definition, the algebra
structure on the vector space $V \ot A$, whose elements
are now denoted by $v \# a$, given by
\[
(v \# a) (u \# b) ~ = ~ v \aa (u) ~\# ~\ab b.
\]
This algebra is defined in such a way that the map
\begin{equation}
\label{eq_t1}
T_1: ~ H = V \# A \lrw End_k(V): ~ T_1(v \# a)(u) ~ = ~ v a(u)
\end{equation}
defines a left representation of $H$ on $V$.

\bigskip
Before describing the Hopf algebroid structure on $H = V \# A$
over $V$, we first give an example of this Theorem.

\begin{exam}
\label{exam_dual}
{\em
The following left action of $D(A)$ on $A^*$ makes $A^*$ into a left
$D(A)$-module algebra:
\begin{equation}
\label{eq_D-onAdual}
D(A) \ni x \ot a: ~ y \Map \xb ~ (a ~\rhu ~y) ~ S^{-1}(\xa), ~~~ y \in A^*.
\end{equation}
Moreover, for this action, Condition (\ref{eq_R}) is satisfied. Indeed,
we have, for any $x, y \in A^*$,
\begin{eqnarray*}
x_t (y) a_t(x) & = & (x_t)_{(2)} y S^{-1}((x_t)_{(1)}) (a_t ~\rhu ~x)\\
& = & x_s y S^{-1}(x_t) ((a_t a_s) ~\rhu ~x) \\
& = & x_s y S^{-1}(x_t) \xa \la \xb, ~ a_t \ra \la x_{(3)}, ~ a_s \ra\\
& = & x_{(3)} y S^{-1}(\xb) \xa \\
& = & xy.
\end{eqnarray*}
Notice that as a Hopf subalgebra of $D(A), ~ A^{*coop}$
acts on $A^*$ from the left by
\begin{equation}
\label{eq_adjoint}
A^{* coop} \ot A^* \lrw A^*: ~~ x \ot y \Map ad_x y ~ =:~ \xb ~ y ~
S^{-1}(\xa).
\end{equation}
It makes $A^*$ into a left $A^{*coop}$-module algebra. We call it
the {\bf left adjoint representation} of $A^{*coop}$ on $A^*$.

As a Hopf subalgebra of $D(A), ~ A$ acts on $A^*$ by
\begin{equation}
\label{eq_left-regular}
A \ot A^* \lrw A^*: ~ a \ot x \Map a \rhu x ~ = ~ \xa \la \xb, ~ a \ra.
\end{equation}
It makes $A^*$ into a left $A$-module algebra.
We call it the {\bf left regular representation} of $A$ on $A^*$.
The corresponding smash product algebra $H = A^* \# A$ is, by definition,
the Heisenberg double of $A^*$. Thus our theorem describes a Hopf algebroid
structure over $A^*$ on its Heisenberg double.
}
\end{exam}

\bigskip
Returning to Theorem \ref{thm_main}, we now describe the Hopf algebroid
structure on $H = V \# A$.

\begin{itemize}
\item
The source map is defined by
\[
\al: ~ V \lrw H: ~ v \Map v \# 1.
\]
It is clearly an algebra homomorphism.

\item
To define the target map $\beta: V \rw H$, we first restrict the action
of $D(A)$ on $V$ to a left action of $A^*$ on $V$. Since the
inclusion
\[
A^{*coop} = (A^{op})^{*} \lrw D(A): ~ x \Map x \ot 1
\]
is a Hopf algebra inclusion,  $V$ becomes
a left $A^{* coop} \cong (A^{op})^*$-module algebra. The
corresponding co-module map, which will
be the target map so we denote it by $\beta$:
\[
\be: ~ V \lrw V \ot A^{op}: ~ \be(v) ~ = ~ x_t(v) \ot a_t
\]
is an algebra homomorphism from $V$ to the direct product
algebra $V \ot A^{op}$.

\end{itemize}

\begin{lem}
\label{lem_beta}
As a map from $V$ to the algebra $H = V \# A$, the map $\beta$
has the following properties:

1) $\al (v) \be(u) ~ = ~ \be(u) \al(v)$ for all $u, v \in V$. In fact, this
is equivalent to Condition (\ref{eq_R});

2) It is an algebra anti-homomorphism from $V$ to $H = V \# A$.
\end{lem}

\noindent
{\bf Proof} 1) For $u, v \in V$, we have, on the one hand,
\[
\al (v) \be(u) ~ = ~ (v \#1) (x_t(u) \# a_t) ~ = ~ v x_t(u) \# a_t,
\]
and on the other hand,
\begin{eqnarray*}
\be(u) \al(v) & = & (x_t(u) \# a_t) ~ (v \# 1) \\
& = & x_t (u) (a_t)_{(1)} (v) \# (a_t)_{(2)} \\
& = & (x_s x_t)(u) ~ a_s (v) \# a_t\\
& = & x_s (x_t(u)) ~ a_s(v) \# a_t.
\end{eqnarray*}
By Condition (\ref{eq_R}), we have
\[
\be(u) \al(v) ~ = ~ v x_t (u) \# a_t ~ = ~ \al(v) \be(u).
\]
It is easy to see that  conversely 1) implies Condition (\ref{eq_R}).

2) Using Condition (\ref{eq_R}) again, we have, for all $u, v \in V$,
\begin{eqnarray*}
\be(u) \be(v) & = & (x_t(u) \# a_t) ~ (x_s(v) \# a_s) \\
& = & x_t(u) (a_t)_{(1)} ~ (x_s(v)) \# (a_t)_{(2)} a_s \\
& = & (x_{\gamma} x_{\xi})(u) ~ a_{\gamma} ~ (x_s(v)) \# a_{\xi} a_s\\
& = & x_s(v) x_{\xi}(u)\# a_{\xi} a_s \\
& = & (x_t)_{(2)} (v) (x_t)_{(1)} (u) \# a_t\\
& = & x_t (vu) \# a_t\\
& = & \be(vu).
\end{eqnarray*}
Thus $\be$ is an algebra anti-homomorphism from $V$ to $H = V \# A$.
\qed

\begin{itemize}

\item The coproduct: The co-product should be a map from $H$ to
$H \otm H$. In our case, the latter is the quotient space of
$H \ot H$ by the subspace spanned by the following
subset:
\[
\{\be(v) (v_1 \# a_1) \ot (v_2 \# a_2) - (v_1 \# a_1) (v v_2 \# a_2): ~
v, v_1, v_2 \in V, ~ a_1, a_2 \in A \}.
\]
Thus each element  in $H \otm H$ is uniquely represented by an element of
the form $\sum (v_i \# a_i) \ot (1 \# b_i))$. In other words,
we have a map
\begin{equation}
\label{eq_gamma}
\gamma:~ H \otm H \lrw H \ot H:~ (v_1 \# a_1) \ot (v_2 \# a_2) \Map
\be(v_2) (v_1 \# a_1) \ot (1 \# a_2)
\end{equation}
which is a section for the natural projection from $H \ot H$ to $H \otm H$.

Define $\Delta: H \lrw H \otm H$ by
\begin{equation}
\label{eq_Delta-hopf}
\Delta (v \# a) ~ = ~ (v \# \aa) \ot (1 \# \ab).
\end{equation}
The associativity of the co-product for $A$ implies that
$\Delta$ is co-associative. It is also easy to see that $\Delta$
is a $(V, V)$-bimodule map.  It remains to prove the compatibility
of the co-product and the product of $H$.

We first observe that
\[
\gamma(H \otm H) ~ = ~ \{h \ot (1 \# b): ~~ h \in H, ~ b \in A\}
\]
is a subalgebra of $H \ot H$. Moreover, the map
\[
\gamma \circ \Delta: ~ H \lrw \gamma (H \otm H)
\]
is an algebra homomorphism. Thus the following defines
a left $H$-module structure on $\gamma(H \otm H)$:
\[
(v \# a) \cdot (h \ot 1 \# b) ~ = ~ \Delta(v \# a) (h \ot 1 \# b)
{}~ = ~ (v \# \aa ) h \ot (1 \# \ab b).
\]
On the other hand, by identifying $\gamma(H \otm H)$ with
$H \otm H$, the right $(H \ot H)$-module structure on
$H \otm H$ now becomes the following
right $(H \ot H)$-module structure on $\gamma(H \otm H)$:
\[
(h \ot 1 \# b) \cdot (h_1 \ot (v_2 \# a_2))
{}~ = ~ \beta(\ab(v_2)) h h_1 \ot (1 \# b_{(2)} \ab).
\]
These two module structures on $\gamma(H \otm H)$ commute. Indeed,
let $v \# a, ~ h_1, ~ h_2 = v_2 \# a_2 \in H$. We have
\beqa
\left( (v \# a) \cdot (h \ot 1 \# b) \right) \cdot (h_1 \ot h_2) & = &
\beta \left( \ab(\ba (v_2)) \right) (v \# \aa) h h_1 ~ \ot ~ (1 \# \ac \bt a_2)
\\
& = & v \beta \left(\ab(\ba (v_2)) \right) (1 \# \aa) h h_1 ~ \ot ~
(1 \# \ac \bt a_2) \\
& = & v (1 \# \ab) \beta(\ba (v_2)) h h_1 ~ \ot ~ (1 \# \ac \bt a_2) \\
& = & (v \# a) \cdot \left( (h \ot 1 \# b) \cdot (h_1 \ot h_2) \right).
\eeqa
We used  1) in Lemma \ref{lem_beta} and the following fact in
deriving the above identities:

\begin{lem}
\label{lem_beta-1}
For any $u \in V$ and $a \in A$,
\begin{equation}
\label{eq_bu}
\beta(\ab(u)) (1 \# \aa) ~ = ~ (1 \# a) \beta (u).
\end{equation}
\end{lem}

\noindent
{\bf Proof of Lemma \ref{lem_beta-1}.} By the definition of the
map $\beta$, we have
\[
\beta(\ab(u)) (1 \# \aa)  ~ = ~ x_s (\ab(u)) ~ \# ~ a_s \aa \in V \ot A
\]
and
\[
(1 \# a) \beta (u) ~ = ~ \aa (x_s(u)) ~ \# ~ \ab a_s \in V \ot A.
\]
Regarding both elements as in $Hom_k(A^*, V)$, we see that
(\ref{eq_bu}) is equivalent to
\[
\xa (\ab(u)) \la \aa, ~ \xb \ra ~ = ~ \aa (\xb(u))
\la \ab, ~ xa \ra
\]
for all $x \in A^*$. But this follows from the following identity
in $D(A)$ which can be proved directly using the definition of
the product of $D(A)$:
\begin{equation}
\label{eq_xa}
\xa ~ \ab ~ \la \aa, ~ \xb \ra ~ = ~ \aa \xb \la \ab, ~ \ba \ra, ~~~~
\forall a \in A, ~ x \in A^*.
\end{equation}
\qed

We now return to the proof of the compatibility of the
co-product $\delta$ and the product of $A$. We need to show that
the kernel of the map
\[
\Phi: ~ H \ot H \ot H \lrw H \otm H: ~~ h \ot h_1 \ot h_2 \Map (\Delta h)
(h_1 \ot h_2)
\]
is a left ideal in $H \ot H^{op} \ot H^{op}$. Let
$h \ot h_1 \ot h_2 \in H \ot H \ot H$ and let
$\sum_i h_i \ot h_{i}^{'} \ot h_{i}^{''} \in \ker \Phi$.
Then
\beqa
\Phi(\sum_i h h_i \ot h_{i}^{'} h_1 \ot h_{i}^{''} h_2) & = &
\sum_i (\gamma \circ \Delta) (h h_i) \cdot (h_{i}^{'} h_1 \ot h_{i}^{''} h_2)\\
& = & \sum_i \left( \Delta(h) \cdot (\gamma \circ \Delta)(h_i) \right)
\cdot (h_{i}^{'} \ot h_{i}^{''}) (h_1 \ot h_2) \\
& = &  \sum_i \Delta(h) \cdot \left((\gamma \circ \Delta)(h_i)
\cdot (h_{i}^{'} \ot h_{i}^{''}) (h_1 \ot h_2) \right) \\
& = & \Delta(h) \cdot \left( \Phi(\sum_i h_i \ot h_{i}^{'} \ot h_{i}^{''})
\cdot (h_1 \ot h_2) \right)\\
& = & 0.
\eeqa
Here we have just used the fact that the left $H$-module
structure and the right $H \ot H$-module
structure on $\gamma(H \otm H)$ commute. Hence $\ker \Phi$ is
a left ideal of $H \ot H^{op} \ot H^{op}$.

\item
The co-unit: define
\begin{equation}
\label{eq_epsilon}
\epsilon:~ H \lrw V: ~ \epsilon(v \# a) ~ = ~ \epsilon(a) v.
\end{equation}
Then it satisfies all the requirements for a co-unit. For example,
the kernel of $\epsilon$ is the space of annihilators in $H$ of the
identity element $1_{\scriptscriptstyle V}$ of $V$ under the left action
$T_1$ of $H$ on $V$. Therefore it is a left ideal of $H$.

\bigskip
We have thus obtained a bi-algebroid structure on the smash-product
algebra $H = V \# A$ over the algebra $V$. The associated
left action of $H$ on $V$ is the one given by (\ref{eq_t1}).

\item
The antipode: introduce a special element $d_0 \in D(A)$ given by
\begin{equation}
\label{eq_d}
d_0 ~ = ~ S^2(a_t) x_t \in D(A).
\end{equation}
We now define the antipode map $\tau$ by
\begin{equation}
\label{eq_tau-main}
\tau: ~ H \lrw H: ~ \tau(v \# a) ~ = ~ (1 \# S(a)) ~ \be (d_0 (v)).
\end{equation}

\end{itemize}

\begin{lem}
\label{lem_d0}
The element $d_0 \in D(A)$ has the following properties:

1) It is invertible in $D(A)$ and $d^{-1}_{0} = S^{-1}(a_t) x_t$.

2) For any $x \ot a \in D(A)$ we have
\begin{equation}
\label{eq_s-sqr}
S_{\scriptscriptstyle D(A)}^{2} (x \ot a) ~ = ~ d_0 (x \ot a) d_{0}^{-1},
\end{equation}
where $S_{\scriptscriptstyle D(A)}$ is the antipode map for $D(A)$
given by $S_{\scriptscriptstyle D(A)}(x \ot a) = S(a) S^{-1}(x)$.

3) $\Delta d_0  =  (R^{21} R) (d_0 \ot d_0)$, where $R$ is the $R$-matrix
for $D(A)$ given by (\ref{eq_R-dfn}) and $R^{21}$ is its flip, i.e.,
\[
R^{21} ~ = ~ (x_t \ot 1 ) ~ \ot ~ (1 \ot a_t) \in D(A) \ot D(A).
\]

4) $d_0$ acts as an algebra isomorphism on $V$.
\end{lem}

\noindent
{\bf Proof.} The first three properties of $d_0$ were proved in
\cite{mt:book} \cite{dr:almost}. We only prove 4).
Condition (\ref{eq_R}) implies,
on the one hand, $m^{op}_{\scriptscriptstyle V} = m_{\scriptscriptstyle V}
\circ R^{-1}$, and on the other hand, $m^{op}_{\scriptscriptstyle V}
 = m_{\scriptscriptstyle V} \circ
R^{21}$. Therefore, we have
\[
m_{\scriptscriptstyle V}  ~ = ~ m_{\scriptscriptstyle V} \circ (R^{21}R).
\]
Hence, using the ``product rule" for the action of $D(A)$ on $V$, we get,
for any $u, v \in V$,
\begin{eqnarray*}
d_0 (u v) & = & m_{\scriptscriptstyle V} \Delta d_0 (u \ot v) \\
& = & m_{\scriptscriptstyle V} \circ  R^{21}R (d_0(u) \ot d_0(v)) \\
& = & d_0(u) d_0(v).
\end{eqnarray*}
Thus $d_0$ acts as an isomorphism on $V$.
\qed

\begin{lem}
\label{lem_antipode-main}
1) $\tau$ is an algebra anti-isomorphism for $H$, and $\tau^{-1}: H \rw H$
is given by
\begin{equation}
\label{eq_tau-inverse}
\tau^{-1} (v \# a) ~ = ~ (1 \# S^{-1}(a)) ~ \be (v) ~ = ~ (1 \# S^{-1}(a))
{}~ (x_t (v) \# a_t).
\end{equation}

2) $\tau \be = \al$ and $\tau \al = \be d_0$.

3) $m_{\scriptscriptstyle H} (\tau \ot id) \Delta = \be \epsilon \tau$.

4) $m_{\scriptscriptstyle H} (id \ot \tau) \gamma \Delta = \al \epsilon$,
where $\gamma: H \otm H \rw H \ot H$ is given by
(\ref{eq_gamma}).
\end{lem}

\noindent
{\bf Proof.} 1) To prove that $\tau$ is an algebra anti-isomorphism, it
is sufficient to prove the following: for any $v \in V$ and $a \in A$:
\[
\tau((1 \#a) (v \#1)) ~ = ~ \tau(v \#1) ~ \tau (1 \# a).
\]
Now
\begin{eqnarray*}
lhs & = & \tau(\aa (v) \# \ab) \\
& = & (1 \# S(\ab)) ~ \be (d_0 \aa(v)) \\
& = & \be(S(\ab) d_0 \aa (v) ) ~ S(\ac).
\end{eqnarray*}
Using (\ref{eq_s-sqr}), we have
\[
S(\ab) d_0 \aa ~ = ~ S(\ab) S^2(\aa) d_0 ~ = ~ \epsilon(\aa) d_0.
\]
Thus
\[
lhs ~ = ~ \be(d_0 (v)) ~ S(a) ~ = ~ rhs.
\]

The proof of 3) uses Identity (\ref{eq_bu}). The proof of
$\tau \alpha = \beta d_0$ is trivial, and so is the proof of 4).
We now give the proof of $\tau \beta = \alpha$. Let $v \in V$ be
arbitrary. Then
\beqa
\tau \beta(v) & = & \tau(x_t(v) ~ \# a_t) \\
& = & (1 ~ \# ~ S(a_t)) ~ \beta (d_0 x_t (v))\\
& = & (1 ~ \# ~ S(a_t)) ~ \left( (x_s d_0 x_t)(v) ~ \# ~ a_s \right)\\
& = & \left( S(a_{\eta}) x_s d_0 x_{\xi} x_{\eta} \right)(v)
{}~ \# ~ S(a_{\xi}) a_s.
\eeqa
By (\ref{eq_s-sqr}),
\[
d_0 x_{\xi} ~ = ~ S^{-2}(x_{\xi}) d_0.
\]
Moreover,
\[
x_s S^{-2}(x_{\xi}) ~ \ot ~ S(a_{\xi}) a_s ~ = ~ 1 \ot 1 \in A^* \ot A.
\]
This can be seen by pairing both sides with arbitrary
elements $a \ot x \in A \ot A^*$. Thus
\beqa
\tau \beta (v) & = & (s(a_{\eta}) d_0 x_{\eta}) (v) ~ \# ~ 1 \\
& = & (S(a_{\eta}) S^{-2}(x_{\eta}) d_0)(v) ~ \# ~ 1.
\eeqa
But
\[
S(a_{\eta}) S^{-2}(x_{\eta}) ~ = ~ S_{\scriptscriptstyle D(A)}^{2} (
S^{-1}(a_{\eta}) x_{\eta}) ~ = ~
S_{\scriptscriptstyle D(A)}^{2} (d_{0}^{-1}) ~ = ~ d_{0}^{-1}.
\]
Hence $\tau \beta (v) = v \# 1 = \alpha(v)$.
\qed

We have finished the proof of Theorem \ref{thm_main}.

\begin{cor}
\label{cor_theta}
The algebra automorphism $\theta$ on $V$ as defined in Proposition
\ref{prop_d} is given by the action of $d_0$ on $V$.
\end{cor}

The proof
of the following Proposition is straightforward.

\begin{prop}
\label{prop_main-T}
Let $V, A$ and $H = V \# A$ be as in Theorem
\ref{thm_main}.
The left representation of $H$ on $V$ associated to the
Hopf algebroid structure on $H$ (see Definition \ref{dfn_induced})
is the one given by (\ref{eq_t1}).
\end{prop}

\section{The Hopf algebroid $H(A^*)$}
\label{sec_hsberg}

In this section, we study in more details the Hopf algebroid
$H(A^*)$ for any finite dimensional Hopf algebra $A$.

\bigskip
The representation $T_1$ of $H(A^*)$ on $A^*$ defined by
(\ref{eq_t1}) now becomes
\[
T_1: ~ H(A^*): ~ \lrw End_k(A^*): ~~
T_1(x \# a) (y ) ~ = ~ x (a ~ \rhu ~ y).
\]
The following Proposition says that $T_1$ is an algebra
isomorphism from $H(A^*)$ to $End_k(A^*)$.

\begin{prop}
\label{prop_iso}
For any $\phi \in End_k (A^*)$,
\begin{equation}
\label{eq_ht1}
T_1 \left(\phi(x_s) x_t ~ \# S^{-1}(a_t) a_s\right) ~ = ~ \phi,
\end{equation}
where $\{a_t\}$ is a basis for $A$ and $\{x_t\}$ is its dual basis for $A^*$.
\end{prop}

\noindent
{\bf Proof.} For any $x \in A^*$,
\beqa
T_1\left(\phi(x_s) x_t ~ \# S^{-1}(a_t) a_s\right)(x) & = & \phi (x_s) x_t
\left(
S^{-1} (a_t) a_s ~ \rhu ~ x\right) \\
& = & \phi (x_s) x_t \xa \la S^{-1}(a_t), ~ \xb \ra \la a_s, ~ \xc \ra\\
& = & \phi(\xc) S^{-1} (\xb) \xa\\
& = & \phi(x).
\eeqa
\qed

\begin{cor}
\label{cor_D-to-H}
The map
\begin{equation}
\label{eq_D-to-H}
D(A) \lrw H(A^*): ~~ a ~ \# ~ x \Map
\xb (a \rhu x_t) S^{-1} (\xa) x_s ~ \# ~ S^{-1}(a_s) a_t
\end{equation}
is an algebra homomorphism.
\end{cor}

\bigskip
The structure maps for the Hopf algebroid $H(A^*)$ are now given by

1) $~\alpha: ~ A^* \rightarrow H(A^*): ~~ x \mapsto x \# 1;$

2) $~ \beta: ~ A^* \rightarrow H(A^*): ~~ x \mapsto x_s x x_t ~ \# S^{-1}(a_t)
a_s;$

3) $~ \epsilon: ~ H(A^*) \rightarrow A^*: ~~ x \# a \mapsto \epsilon(a) x;$

4) $~ \Delta: ~ H(A^*) \rightarrow H(A^*) \ot_{A^*} H(A^*): ~~
x \# a \mapsto x \# \aa ~ \ot ~ 1 \# \ab;$

5) $~ \tau: ~ H(A^*) \rightarrow H(A^*): ~ x \# a \mapsto
\left(S(a) \rhu (x_s S^{-2}(x)) \right) x_t ~ \# ~ S^{-1}(a_t) a_s.$

\section{An example of a quantum groupoid associated to $SL_q(2)$}
\label{sec_q-poid}

Let $d >1$ be an odd integer. Let $q$ be a root of unity of order $d$.
Let $A$ be the Hopf algebra with two generators $E$ and $K$ and
relations:
\beqa
& & KE ~ = ~ q^2 EK, ~~ K^d ~ = ~ 1, ~~ E^d ~ = ~ 0;\\
& & \Delta K ~ = ~ K \ot K, ~~ \Delta E ~ = ~ 1 \ot E ~ + ~ E \ot K; \\
& & S(K) ~ = ~ K^{-1}, ~~ S(E) ~ = ~ -EK^{-1};\\
& & \epsilon (K) ~ = ~ 1, ~~ \epsilon (E) ~ = ~ 0.
\eeqa
It is easy to see that the set $\{E^m K^n \}_{m, n = 0, 1, ..., d-1}$
is a basis for $A$, so $\dim A = d^2$.

The dual Hopf algebra $A^*$ of $A$ is the Hopf algebra with
two generators $\kappa$ and $\eta$ and relations
\beqa
& & \kappa \eta ~ = ~ q^{-2} \eta \kappa, ~~ \kappa^d ~ = ~ 1, ~~
\eta^d ~ = ~ 0;\\
& & \Delta \kappa ~ = ~ \kappa \ot \kappa, ~~ \Delta \eta ~ = ~
\eta \ot 1 ~ + ~ \kappa \ot \eta;\\
& & S(\kappa) ~ = ~ {\kappa}^{-1}, ~~ S(\eta) ~ = ~ -\kappa^{-1}
\eta;\\
& & \epsilon (\kappa) ~ = ~ 1, ~~ \epsilon (\eta) ~ = ~ 0.
\eeqa
Again, the set $\{ \eta^i \kappa^j\}_{i,j = 0, 1, ..., d-1}$ forms
a basis
for $A^*$, with the pairing between $A$ and $A^*$ given by
\begin{equation}
\label{eq_pairing}
\la \eta^i \kappa^j, ~ E^m K^n \ra ~ = ~ \delta_{mi} (i)_{q^2} ! q^{2j(i+n)},
\end{equation}
where
\[
(i)_{q^2} ~ = ~ {\frac{q^{2i} - 1}{q^2 -1}}
\]
and
\[
(i)_{q^2} ! ~ = ~ (i)_{q^2} (i-1)_{q^2} \cdots (2)_{q^2} (1)_{q^2}.
\]
By definition,
\[
(0)_{q^2} ~ = ~ 0, ~~~~~~ (0)_{q^2} ! ~ = ~ 1.
\]
In particular,
\beqa
& & \la \kappa, ~ E^m \ra ~ = ~ \delta_{m0}, ~~~~ \la \kappa, ~ K^n \ra ~ = ~
q^{2n},\\
& & \la \eta, ~ E^m \ra ~ = ~ \delta_{m1}, ~~~~ \la \eta, ~ K^n \ra ~ = ~ 0.
\eeqa

The
Drinfeld Double $D(A)$ of $A$ has four generators $E, K, \eta$ and
$\kappa$ with the following relations:
\beqa
& & K^d ~ = ~ 1, ~~~~ \kappa^d ~ = ~ 1, ~~~~ E^d ~ = ~ 0, ~~~~ \eta^d ~ = ~
0;\\
& & K E ~ = ~ q^2 EK, ~~~ ~\kappa \eta ~ = ~ q^{-2} \eta \kappa;\\
& & K \kappa ~ = ~ \kappa K, ~~~~ E \kappa ~ = ~ q^{-2} \kappa E,\\
& & K \eta ~ = ~ q^{-2} \eta K, ~~~~ E \eta ~ = ~ q^{-2}
(-1 + \eta E + \kappa K);\\
& & S(K) ~ = ~ K^{-1}, ~~~~ S(\kappa) ~ = ~ \kappa^{-1}, ~~~~
S(E) ~ = ~ -E K^{-1}, ~~~~ S(\eta) ~ = ~ - \eta \kappa^{-1};\\
& & \Delta K ~ = ~ K \ot K, ~~~~ \Delta \kappa ~ = ~ \kappa \ot \kappa, \\
& & \Delta E ~ = ~ 1 \ot E ~ + ~ E \ot K, ~~~~ \Delta \eta ~ = ~ 1 \ot \eta ~ +
{}~
\eta \ot \kappa;\\
& & \epsilon(K) ~ = ~ \epsilon (\kappa) ~ = ~ 1, ~~~~ \epsilon(E) ~ = ~
\epsilon (\eta) ~ = ~ 0.
\eeqa

The $R$-matrix for $D(A)$ is given by
\begin{equation}
\label{eq_R-1}
R ~ = ~ {\frac{1}{d}} \sum_{0 \leq m, n, j \leq d-1} {\frac{1}{(m)_{q^2}!}}
q^{-2j(m+n)} E^m K^n ~ \ot ~ \eta^m \kappa^j \in D(A) \ot D(A).
\end{equation}

The left regular representation of $A$ on $A^*$ (see
Example \ref{exam_dual}) is given by
\beqa
& & K ~ \rhu ~ (\eta^i \kappa^j) ~ = ~ q^{2j} \eta^i \kappa^j \\
& & E ~ \rhu ~ (\eta^i \kappa^j) ~ = ~ (i)_{q^2} q^{2(j-i+1)} \eta^{i-1}
\kappa^{j+1}.
\eeqa
In particular,
\beqa
& & K ~ \rhu ~ \kappa ~ = ~ q^2 \kappa, ~~~~
K ~ \rhu ~ \eta ~ = ~ \eta; \\
& & E ~ \rhu ~ \kappa ~ = ~ 0, ~~~~
E ~ \rhu ~ \eta ~ = ~ \kappa.
\eeqa
This makes $A^*$ into a left $A$-module algebra.

The left adjoint representation of $A^{* coop}$ on $A^*$
(see Example \ref{exam_dual}) is given by
\beqa
& & ad_{\kappa} (\eta^i \kappa^j) ~ = ~ q^{-2i} \eta^i \kappa^j \\
& & ad_{\eta} (\eta^i \kappa^j) ~ = ~ (1-q^{-2j}) \eta^{i+1}
\kappa^{j-1}.
\eeqa
In particular,
\beqa
& & ad_{\kappa} \kappa ~ = ~ \kappa, ~~~~
ad_{\kappa} \eta ~ = ~ q^{-2} \eta; \\
& & ad_{\eta} \kappa ~ = ~ (1 - q^{-2}) \eta, ~~~~
ad_{\eta} \eta ~ = ~ 0.
\eeqa
This action makes $A^*$ into a left $A^{*coop}$-module
algebra.

\bigskip
The left regular representation of $A$ on $A^*$ and the
left adjoint representation of $A^{*coop}$ on $A^*$ combine to
give a left action of $D(A)$ on $A^*$ (see Example
\ref{exam_dual}). This makes $A^*$ into a left
$D(A)$-module algebra. Thus by Theorem \ref{thm_main},
there is a Hopf algebroid structure over $A^*$  on the Heisenberg
double of $A^*$.

\bigskip
The Heisenberg double $H(A^*) = A^* \# A$ of $A^*$
is the algebra with four generators, again denoted by $E, K, \kappa$ and
$\eta$ with relations:
\beqa
& & KE ~ = ~ q^2 EK, ~~~ K^d ~ = ~ 1, ~~~ \kappa^d ~ = ~ 0;\\
& & \kappa \eta ~ = ~ q^{-2} \eta \kappa, ~~~ \kappa^d ~ = ~ 1, ~~~
\eta^d ~ = ~ 0;\\
& & K \kappa ~ = ~ q^2 \kappa K, ~~~ K \eta ~ = ~ \eta K, ~~~ E \kappa ~ = ~
\kappa E, ~~~ E \eta ~ = ~ \eta E ~ + ~ \kappa K.
\eeqa
It has a natural left representation $T_1$ on $A^*$ given by
(\ref{eq_t1}). Namely, under $T_1, ~ A$ acts on $A^*$ by the
left regular representatin and $A^*$ acts on $A^*$ by left
multiplications.

\bigskip
The Hopf algebroid structure on $H(A^*)$ over $A^*$ is the following:

1) The source map
\beqa
\alpha: ~ A^* \lrw H(A^*) ~ =~ A^* ~ \# ~ A: & ~~& \alpha(k) ~ = ~ k ~ \#  1\\
& ~~ & \alpha(\eta) ~ = ~ \eta ~ \# ~ 1;
\eeqa

2) The target map
\beqa
\beta: ~ A^* \lrw H(A^*) ~ = ~ A^* ~ \# ~ A: & ~~ &
\beta(k) ~ = ~ k  ~ + ~ (1-q^{-2}) \eta E K^{-1}\\
& ~~ & \beta(\eta) ~ = ~ \eta K^{-1}.
\eeqa

3) Identify the space $H \otm H$ with the space spanned by the set
\[
\{ \eta^i \kappa^j E^m K^n ~ \ot ~ E^r K^s: ~~ i, j, m, n, r, s = 0, ...,
d-1\},
\]
then the co-product on $H(A^*)$ is given by
\[
\Delta: ~ H(A^*) \lrw H \otm H: ~~~
\eta^i \kappa^j E^m K^n \Map \sum_r \left(\begin{array}{c}m\\r \end{array}
\right)_{q^2} \eta^i \kappa^j E^m K^n ~
\ot ~ E^r K^{m+n-r}.
\]

4) The antipode map
\beqa
\tau: ~ H(A^*) \lrw H(A^*): & ~~& \tau(k) ~ = ~ q^2 k ~ + ~ (q^2 - 1) \eta E
K^{-1}\\
& ~~ & \tau(\eta) ~ =~ \eta K^{-1} \\
& ~~& \tau(K) ~ = ~ K^{-1} \\
& ~~& \tau(E) ~ = ~ -E K^{-1}.
\eeqa

5) The co-unit map
\[
\epsilon: ~ H(A^*) ~ = ~ A^* ~ \# ~ A \lrw A^*: ~~ \epsilon(\eta^i \kappa^j
E^m K^n) ~ = ~ \delta_{m0} \eta^i \kappa^j.
\]

6) The induced algebra automorphism $\theta$ on $A^*$ as defined
in Proposition \ref{prop_d} is given by
\beqa
\theta ~ = ~ q^2 S^{-2}: & ~~& \theta(k) ~ = ~ q^2 k\\
& ~~ & \theta(\eta) ~ = ~ \eta.
\eeqa

\end{document}